\documentclass[a4paper]{amsart}

\usepackage{epsfig}
\usepackage{amsthm,amsfonts}
\usepackage{amssymb,graphicx,color}
\usepackage[all]{xy}
\usepackage{cancel} %\cancel{}, \bcancel{}, \xcancel{}
\usepackage{verbatim}
\usepackage{hyperref}
\usepackage{enumitem}
\usepackage{tikz-cd}

\newtheorem{innercustomcor}{Corollary}
\newenvironment{customcor}[1]
  {\renewcommand\theinnercustomcor{#1}\innercustomcor}
  {\endinnercustomcor}
\newtheorem{innercustomprop}{Proposition}
\newenvironment{customprop}[1]
  {\renewcommand\theinnercustomprop{#1}\innercustomprop}
  {\endinnercustomprop}
\newtheorem{theorem}{Theorem}[section]
\newtheorem*{theorem*}{Theorem}

\newtheorem{corollary}[theorem]{Corollary}
\newtheorem{proposition}[theorem]{Proposition}
\newtheorem{claim}{Claim}[theorem]

\theoremstyle{definition}
\newtheorem{definition}[theorem]{Definition}

\newtheorem{conjecture}{Conjecture}
\newtheorem*{conjecturezmodd}{Conjecture Z mod 2}

\theoremstyle{remark}

\newcommand{\R}{\mathbb{R}}
\newcommand{\C}{\mathbb{C}}

\begin{document}

\title[On the Milnor fibres of initial forms]
{On the Milnor fibres of initial forms of topologically equivalent holomorphic functions}

\author[J. Edson Sampaio]{Jos\'e Edson Sampaio}
\address{Jos\'e Edson Sampaio:  
              Departamento de Matem\'atica, Universidade Federal do Cear\'a,
	      Rua Campus do Pici, s/n, Bloco 914, Pici, 60440-900, 
	      Fortaleza-CE, Brazil. \newline   
              E-mail: {\tt edsonsampaio@mat.ufc.br}                    
}

\thanks{The author was partially supported by CNPq-Brazil grant 310438/2021-7. This work was supported by the Serrapilheira Institute (grant number Serra -- R-2110-39576). 
}
%\thanks{The authors were partially supported by CNPq-Brazil}

\begin{abstract}
    Budur, Fernandes de Bobadilla, Le and Nguyen (2022) conjectured that if two germs of holomorphic functions are topologically equivalent, then the Milnor fibres of their initial forms are homotopy equivalent. In this note, we give affirmative answers to this conjecture in the case of plane curves. We show also that a positive answer to this conjecture implies in a positive answer to the famous Zariski multiplicity conjecture both in the case of right equivalence or in the case of hypersurfaces with isolated singularities. 
\end{abstract}

\maketitle
\tableofcontents

\section{Introduction}

Let $f\colon (\mathbb{C}^n,0)\to (\mathbb{C},0)$ be the germ of a holomorphic function at the origin and we write
$$f=f_m+f_{m+1}+\cdots+f_k+\cdots,$$ 
where each $f_k$ is identically zero or is a homogeneous polynomial of degree $k$ and $f_m\neq 0$. In this case, $m={\rm ord}_0 f$ and the polynomial $f_m$ is called {\bf the initial form of $f$} and denoted by $f_*$. Let $(V(f),0)$ be the germ of the zero set of $f$ at origin. When $f$ is a reduced germ, we recall that \emph{the multiplicity of} $V(f)$ at the origin, denoted by $m(V(f),0)$, is defined as $m(V(f),0):= m.$

The famous Zariski multiplicity conjecture asserts that the multiplicity of a hypersurface is a topological invariant. More precisely, we have the following:
\begin{conjecture}\label{conj:zariski}[Zariski multiplicity conjecture]
Let $f,g\colon(\mathbb{C}^n,0)\to (\mathbb{C},0)$ be two reduced holomorphic function-germs. If there is a homeomorphism $\varphi\colon(\mathbb{C}^n,V(f),0)\to (\mathbb{C}^n,V(g),0)$, then $m(V(f),0)=m(V(g),0)$.
\end{conjecture}

The above conjecture was proposed by Zariski as a question in \cite{Zariski:1971}. More than 50 years later, this is still an open problem even in the case of germ of complex analytic functions with isolated singularities. Recently, Zariski's Multiplicity Conjecture for families with isolated singularities was confirmed by Fern\'andez de Bobadilla and Pe\l ka \cite{BobadillaP:2022}.

Here, we also consider the following weaker versions of Conjecture \ref{conj:zariski}:

\begin{conjecture}\label{conj:zariski_right}
Let $f,g\colon(\mathbb{C}^n,0)\to (\mathbb{C},0)$ be two reduced holomorphic function-germs. If there is a homeomorphism $\varphi\colon(\mathbb{C}^n,0)\to (\mathbb{C}^n,0)$ such that $f=g\circ \varphi$, then $m(V(f),0)=m(V(g),0)$.
\end{conjecture}

\begin{conjecture}\label{conj:zariski_isolated}
Let $f,g\colon(\mathbb{C}^n,0)\to (\mathbb{C},0)$ be two holomorphic function-germs with isolated singularity (at the origin). If there is a homeomorphism $\varphi\colon(\mathbb{C}^n,V(f),0)\to (\mathbb{C}^n,V(g),0)$, then $m(V(f),0)=m(V(g),0)$.
\end{conjecture}

Recently, in \cite{BudurFLN:2020} the following conjecture was proposed:

\begin{conjecture}\label{conj:Milnor_fiber_initial_form}
Let $f,g\colon(\mathbb{C}^n,0)\to (\mathbb{C},0)$ be two holomorphic function-germs. If there is a homeomorphism $\varphi\colon(\mathbb{C}^n,0)\to (\mathbb{C}^n,0)$ such that $f=g\circ \varphi$, then $f_*^{-1}(1)$ and $g_*^{-1}(1)$ are homotopy equivalent.
\end{conjecture}

Here, we also consider the following weaker versions of Conjecture \ref{conj:Milnor_fiber_initial_form}:

\begin{conjecture}\label{conj:Milnor_fiber_initial_form_homology}
Let $f,g\colon(\mathbb{C}^n,0)\to (\mathbb{C},0)$ be two holomorphic function-germs. If there is a homeomorphism $\varphi\colon(\mathbb{C}^n,0)\to (\mathbb{C}^n,0)$ such that $f=g\circ \varphi$, then $f_*^{-1}(1)$ and $g_*^{-1}(1)$ have the same homology.
\end{conjecture}

\begin{conjecture}\label{conj:Milnor_fiber_initial_form_isolated}
Let $f,g\colon(\mathbb{C}^n,0)\to (\mathbb{C},0)$ be two holomorphic function-germs with isolated singularity (at the origin). If there is a homeomorphism $\varphi\colon(\mathbb{C}^n,0)\to (\mathbb{C}^n,0)$ such that $f=g\circ \varphi$, then $f_*^{-1}(1)$ and $g_*^{-1}(1)$  have the same homology.
\end{conjecture}

Following the ideas in \cite{Sampaio:2020a}, we obtain the following result, which presents an interesting relation between Conjectures \ref{conj:Milnor_fiber_initial_form} and \ref{conj:zariski}. 
As was already mentioned above, in \cite{BobadillaP:2022} it was proved that Conjecture \ref{conj:zariski_isolated} 
has a positive answer in the case of families and with that it was proved that Conjecture \ref{conj:Milnor_fiber_initial_form_isolated} also has a positive answer (see \cite[Remark 8.2]{BobadillaP:2022}). Here, we show that the converse is true as well even in the case of pair of functions.
\begin{customprop}{\ref*{prop:conj_initial_implies_zariski}}
If Conjecture \ref{conj:Milnor_fiber_initial_form_homology} has a positive answer, then Conjecture \ref{conj:zariski_right} also has a positive answer.
\end{customprop}
\begin{customprop}{\ref*{prop:conj_initial_implies_zariski_isolated}}
If Conjecture \ref{conj:Milnor_fiber_initial_form_isolated} has a positive answer, then Conjecture \ref{conj:zariski_isolated} also has a positive answer.
\end{customprop}

% In this case, the tangent cokne of $V(f)$ at the origin is $C(V(f),0)=V(f_*)$.
The above proposition says, in particular, that Conjecture \ref{conj:Milnor_fiber_initial_form} is possibly harder than Conjecture \ref{conj:zariski}. 
However, the next result shows that, like Conjecture \ref{conj:zariski}, Conjecture \ref{conj:Milnor_fiber_initial_form} also has a positive answer in dimension 2.

\begin{customprop}{\ref*{prop:conj_initial_dim_two}}
 Let $f,g\colon (\C^2,0)\to (\C,0)$ be two (not necessarily reduced) holomorphic function-germs. If there is a homeomorphism $\varphi\colon(\mathbb{C}^2,V(f),0)\to (\mathbb{C}^n,V(g),0)$, then the Milnor fibres of $f_*$ and $g_*$ are homotopy equivalent.
\end{customprop}

\begin{definition}
Let $\varphi\colon U\subset \R^n\to \R^k$ be a mapping such that $0\in U$ and $\varphi(0)=0$ and let $T=\{t_j\}_j\subset (0,+\infty)$ be a sequence such that $\lim t_j=0$. We say that $\varphi$ {\bf has a pseudo-derivative at $0$ with respect to $T$}, $\psi\colon \R^n\to \R^k$ with $\psi(0)=0$, if $\varphi$ is continuous at $0$ and there is $r>0$ such that $B_r(0)\subset U$ and the sequence of mappings $\varphi_j\colon B_{\frac{r}{t_j}}(0)\to \R^k$, given by $\varphi_j(v)=\frac{\varphi(t_jv)}{t_j}$, uniformly converges on compacts to $\psi$.
\end{definition}

We also have the following result:

\begin{customprop}{\ref*{prop:weak_bi-Lip_A-equiv}}
 Let $f,g\colon (\C^n,0)\to (\C,0)$ be two (not necessarily reduced) holomorphic function-germs. Let $\varphi\colon (\C^n,0)\to (\C^n,0)$ and $\phi\colon (\C,0)\to (\C,0)$ be two mappings such that $f=\phi\circ g\circ \varphi$. Assume that $\varphi$ (resp. $\phi$) has a pseudo-derivative at $0$ with respect to a sequence $T=\{t_j\}_j\subset (0,+\infty)$ (resp. $T'=\{t_j^k\}_j\subset (0,+\infty)$, where $k={\rm ord}_0 g$), $d\varphi\colon \C^n\to \C^n$ (resp. $d\phi\colon \C\to \C$). If $d\varphi$ and $d\phi$ are homeomorphisms, then the Milnor fibres of $f_*$ and $g_*$ are homotopy equivalent.
\end{customprop}

As a consequence of Proposition \ref{prop:weak_bi-Lip_A-equiv} and \cite[Theorem 3.2]{Sampaio:2016} (see also \cite{Sampaio:2015}), we obtain the following result:
\begin{customcor}{\ref*{cor:bi-Lip_A_equiv}}
 Let $f,g\colon (\C^n,0)\to (\C,0)$ be two holomorphic function-germs. Assume that there are bi-Lipschitz homeomorphisms $\varphi\colon (\C^n,0)\to (\C^n,0)$ and $\phi\colon (\C,0)\to (\C,0)$ such that $f=\phi\circ g\circ \varphi$. Then the Milnor fibres of $f_*$ and $g_*$ are homotopy equivalent.
\end{customcor}

Since any subanalytic bi-Lipschitz homeomorphism is a blow-spherical homeomorphism (see the definition of blow-spherical homeomorphism in Subsection \ref{subsec:bs}), it is natural to ask if Conjecture \ref{conj:Milnor_fiber_initial_form} holds true for the blow-spherical homeomorphisms. In order to learn more about the local properties of the blow-spherical equivalence see \cite{BirbrairFG:2017}, \cite{Sampaio:2015}, \cite{Sampaio:2020a}, \cite{Sampaio:2020b}, \cite{Sampaio:2021b} and \cite{SampaioS:2024}.

In this paper, we also prove that Conjecture \ref{conj:Milnor_fiber_initial_form} holds true for the blow-spherical homeomorphisms.
\begin{customprop}{\ref*{prop:bs_R-equiv}}
 Let $f,g\colon (\C^n,0)\to (\C,0)$ be two reduced holomorphic function-germs. Let $\varphi\colon (\C^n,0)\to (\C^n,0)$ be a blow-spherical homeomorphism such that $f=g\circ \varphi$. Then the Milnor fibres of $f_*$ and $g_*$ are homotopy equivalent.
\end{customprop}

Related to Conjecture \ref{conj:Milnor_fiber_initial_form_isolated}, we have the following weaker statement of the Zariski's Problem B:

\begin{conjecture}\label{conj:weaker_problem_b}
Let $f,g\colon(\mathbb{C}^n,0)\to (\mathbb{C},0)$ be two reduced holomorphic function-germs. If there is a homeomorphism $\varphi\colon(\mathbb{C}^n,0)\to (\mathbb{C}^n,0)$ such that $f=g\circ \varphi$, then $f_*^{-1}(0)\setminus \{0\}$ and $g_*^{-1}(0)\setminus \{0\}$ have the same homology.
\end{conjecture}

\begin{customprop}{\ref*{prop:conj_weaker_problem_b_implies_zariski}}
If Conjecture \ref{conj:weaker_problem_b} has a positive answer, then Conjecture \ref{conj:zariski_right} also has a positive answer.
\end{customprop}

Here, we also consider the following version of the Zariski multiplicity conjecture:

\begin{conjecturezmodd}\label{conj:zariski_mod_d}
Let $f,g\colon(\mathbb{C}^n,0)\to (\mathbb{C},0)$ be two reduced holomorphic function-germs. If there is a homeomorphism $\varphi\colon(\mathbb{C}^n,0)\to (\mathbb{C}^n,0)$ such that $f=g\circ \varphi$, then $m(V(f),0)\equiv m(V(g),0)\, {\rm mod}\, 2$.
\end{conjecturezmodd}

Clearly, if Conjecture \ref{conj:zariski_right} has a positive answer, then Conjecture Z mod 2 also has a positive answer. The following result shows that the reciprocal holds true as well.

\begin{customprop}{\ref*{prop:equiv_zariski}}
If Conjecture Z mod 2 has a positive answer, then Conjecture \ref{conj:zariski_right} also has a positive answer.
\end{customprop}

\section{Some results about Conjecture \ref{conj:Milnor_fiber_initial_form}}

\subsection{Conjecture \ref{conj:Milnor_fiber_initial_form} has a positive answer in dimension 2}
\begin{proposition}\label{prop:conj_initial_dim_two}
 Let $f,g\colon (\C^2,0)\to (\C,0)$ be two (not necessarily reduced) holomorphic function-germs. If there is a homeomorphism $\varphi\colon(\mathbb{C}^2,0)\to (\mathbb{C}^2,0)$ such that $f=g\circ \varphi$, then the Milnor fibres of $f_*$ and $g_*$ are homotopy equivalent.
\end{proposition}
\begin{proof}
Let $f=f_1^{m_1}\cdots f_r^{m_r}$ (resp. $g=g_1^{\tilde m_1}\cdots g_s^{\tilde m_s}$) be the irreducible decomposition of $f$ (resp. $g$). Since there is a homeomorphism $\varphi\colon(\mathbb{C}^2,0)\to (\mathbb{C}^2,0)$ such that $f=g\circ \varphi$, then by \cite[Theorem 0.1]{Parusinski:2008}, $r=s$ and we may assume, possibly after reordering the indices, that $m_i=\tilde m_i$ and $\varphi(V(f_i))=V(g_i)$, for all $i\in\{1,..,r\}$. Let $\tilde f=f_1\cdots f_r$ and $\tilde g=g_1\cdots g_r$. By \cite[Theorem 0.1]{Parusinski:2008}, there is a homeomorphism $\tilde\varphi\colon(\mathbb{C}^2,0)\to (\mathbb{C}^2,0)$ such that $\tilde f=\tilde g\circ \tilde \varphi$ and $\tilde\varphi(V(f_i))=V(g_i)$ for all $i\in\{1,..,r\}$.
By the bi-Lipschitz classification of complex analytic curves, which was completed by Pichon and Neumann \cite{N-P}, with previous contributions by Pham and Teissier
\cite{P-T} (see its published version in \cite{P-T:2020}) and Fernandes \cite{Fernandes:2003}, there is a bi-Lipschitz homeomorphism $\psi\colon(\mathbb{C}^2,0)\to (\mathbb{C}^2,0)$ such that $\psi(V(f_i))=V(g_i)$ for all $i\in\{1,..,r\}$. Let $\ell_1,...,\ell_u\colon \C^2\to \C$ (resp. $\ell_1,...,\ell_v\colon \C^2\to \C$) be linear polynomials such that $C(V(\tilde f),0)=\bigcup \limits_{j=1}^u V(\ell_j)$ (resp. $C(V(\tilde g),0)=\bigcup \limits_{j=1}^v V(\tilde \ell_j)$). After multiplying by a nonzero complex number, we may assume that for each $i\in \{1,...,r\}$ there are nonnegative integer numbers $k_1^i,...,k_u^i$ (resp. $\tilde k_1^i,...,\tilde k_v^i$), not all zero, such that ${f_i}_*=\ell_1^{k_1^i}\cdots \ell_u^{k_u^i}$ (resp. ${g_i}_*=\tilde\ell_1^{\tilde k_1^i}\cdots \tilde \ell_u^{\tilde k_v^i}$).
It follows from \cite[Proposition 1.6]{FernandesS:2016} (see also \cite[Theorem 5.3]{FernandesS:2023}) that $u=v$ and there is a bi-Lipschitz homeomorphism $d\psi\colon  (\C^2,0)\to (\C^2,0)$ such that after reordering the indices, if necessary, $d\psi(V(\ell_j))=V(\tilde\ell_j)$ for each $j\in \{1,...,u\}$ and, moreover, for each $i\in \{1,...,r\}$ we have $k_j^i=\tilde k_j^i$ for all $j\in \{1,...,u\}$.

Since $f_*={f_1}_*^{m_1}\cdots {f_r}_*^{m_r}$ and $\tilde f_*={f_1}_*\cdots {f_r}_*$ and similarly $g_*={g_1}_*^{m_1}\cdots {g_r}_*^{m_r}$ and $\tilde g_*={g_1}_*\cdots {g_r}_*$, by \cite[Theorem 0.1]{Parusinski:2008} again, there is a homeomorphism $\tilde\psi\colon(\mathbb{C}^2,0)\to (\mathbb{C}^2,0)$ such that $f_*=g_*\circ \tilde \psi$ and, in particular, the Milnor fibres of $f_*$ and $g_*$ are homotopy equivalent.
\end{proof}

\subsection{Conjecture \ref{conj:Milnor_fiber_initial_form} implies Conjecture \ref{conj:zariski} for the right equivalence}
Here, we need the following result proved separately by A'Campo in \cite{Acampo:1973} and L\^e in \cite{Le:1973}.
\begin{theorem}[A'Campo-L\^e's Theorem]
Let $f,g\colon(\mathbb{C}^n,0)\to (\mathbb{C},0)$ be two complex analytic functions. Suppose that $V(f)$ is smooth at 0. If there is a homeomorphism $\varphi\colon(\mathbb{C}^n,V(f),0)\to (\mathbb{C}^n,V(g),0)$, then $V(g)$ is also smooth at 0.
\end{theorem}

\begin{proposition}\label{prop:conj_initial_implies_zariski}
If Conjecture \ref{conj:Milnor_fiber_initial_form_homology} (resp. Conjecture \ref{conj:Milnor_fiber_initial_form_isolated}) has a positive answer, then Conjecture \ref{conj:zariski_right} (resp. Conjecture \ref{conj:zariski_isolated}) also has a positive answer.
\end{proposition}
\begin{proof}
Assume that Conjecture \ref{conj:Milnor_fiber_initial_form_homology} has a positive answer.

Suppose by contradiction that Conjecture \ref{conj:zariski_right} does not have a positive answer. Thus, there are two reduced holomorphic function-germs $f,g\colon(\mathbb{C}^n,0)\to (\mathbb{C},0)$ and a homeomorphism $\varphi\colon(\mathbb{C}^n,0)\to (\mathbb{C}^n,0)$ such that $f=g\circ \varphi$, but $m(V(f),0)\not=m((V(g),0)$. Since Conjecture \ref{conj:zariski} has a positive answer when $n=2$, then $n\geq 3$.

We assume that $m:=m(V(f),0)>k:=m((V(g),0)$, the case $m(V(f),0)<m((V(g),0)$ is proved in the same way. 

By A'Campo-L\^e's Theorem, $k>1$. Thus, we define $\tilde f,\tilde g\colon (\mathbb{C}^n\times \C,0)\to (\mathbb{C},0)$ given by $\tilde f(x,t)=f(x)+t^k$ and $\tilde g(x,t)=g(x)+t^k$.  Then, the homeomorphism $\tilde\varphi\colon(\mathbb{C}^n\times \C,0)\to (\mathbb{C}^n\times \C,0)$ given by $\tilde \varphi(x,t)=(\varphi(x),t)$ satisfies $\tilde f=\tilde g\circ \tilde \varphi$. Since we are assuming that Conjecture \ref{conj:Milnor_fiber_initial_form_homology} has a positive answer, $\tilde f_*^{-1}(1)$ and $\tilde g_*^{-1}(1)$ have the same homology. But $\tilde g_*$ is a reduced polynomial, and thus $\tilde g_*^{-1}(1)$ is connected. However, $\tilde f_*(x,t)=t^k$, and thus $\tilde f_*^{-1}(1)$ is not connected, which is a contradiction. Therefore, Conjecture \ref{conj:zariski_right} has a positive answer.
\end{proof}

We also have the following:
\begin{proposition}\label{prop:conj_initial_implies_zariski_isolated}
If Conjecture  \ref{conj:Milnor_fiber_initial_form_isolated} has a positive answer, then Conjecture \ref{conj:zariski_isolated} also has a positive answer.
\end{proposition}
\begin{proof}
Assume that Conjecture \ref{conj:Milnor_fiber_initial_form_isolated} has a positive answer.

Suppose by contradiction that Conjecture \ref{conj:zariski_isolated} does not have a positive answer. Thus, there are two holomorphic function-germs $f,g\colon(\mathbb{C}^n,0)\to (\mathbb{C},0)$ with isolated singularities and a homeomorphism $\psi\colon(\mathbb{C}^n,0)\to (\mathbb{C}^n,0)$ such that $\psi(V(f))=V(g)$, but $m(V(f),0)\not=m((V(g),0)$. 

Let $g_1=c_1\circ g\circ c_n$, where $c_k\colon \C^k\to \C^k$ is given by $c_k(z_1,...,z_k)=(\overline{z_1},...,\overline{z_k})$ and $\overline{z_j}$ denotes the complex conjugation of $z_j$. Note that $m(V(g_1),0)=m(V(g),0)$.
Using the corollary in \cite{King:1978} and Corollary 2 in \cite{Saeki:1989}, there is a homeomorphism $\varphi\colon(\mathbb{C}^n,0)\to (\mathbb{C}^n,0)$ such that $f=g\circ \varphi$ or $f=g_1\circ \varphi$. Thus, changing $g$ by $g_1$, if necessary, we assume that $f=g\circ \varphi$. Now, we can argue as in the proof of Proposition \ref{prop:conj_initial_implies_zariski} to obtain a contradiction.

Therefore, Conjecture \ref{conj:zariski_isolated} also has a positive answer.
\end{proof}

\subsection{Case of a weak version of the bi-Lipschitz right-left equivalence}

\begin{proposition}\label{prop:weak_bi-Lip_A-equiv}
 Let $f,g\colon (\C^n,0)\to (\C,0)$ be two (not necessarily reduced) holomorphic function-germs. Let $\varphi\colon (\C^n,0)\to (\C^n,0)$ and $\phi\colon (\C,0)\to (\C,0)$ be two mappings such that $f=\phi\circ g\circ \varphi$. Assume that $\varphi$ (resp. $\phi$) has a pseudo-derivative at $0$ with respect to a sequence $T=\{t_j\}_j\subset (0,+\infty)$ (resp. $T'=\{t_j^k\}_j\subset (0,+\infty)$, where $k={\rm ord}_0 g$), $d\varphi\colon \C^n\to \C^n$ (resp. $d\phi\colon \C\to \C$). If $d\varphi$ and $d\phi$ are homeomorphisms, then the Milnor fibres of $f_*$ and $g_*$ are homotopy equivalent.
\end{proposition}
\begin{proof}
We claim $f_*=d\phi\circ g_*\circ d\varphi$. Indeed, since we have uniform convergence on compacts in the definition of the pseudo-derivatives and $d\phi$ is a continuous function, we obtain that
$$
\lim\limits_{j\to +\infty}\frac{\phi\circ g\circ \varphi(t_jv)}{t_j^k}=\lim\limits_{j\to +\infty}\frac{\phi\left(t_j^k\frac{g\circ \varphi(t_jv)}{t_j^k}\right)}{t_j^k}=d\phi\circ g_*\circ d\varphi(v).
$$
for all $v\in \C^n$. In particular, the above limit is finite for any $v\in \C^n$. This implies ${\rm ord}_0f\geq k={\rm ord}_0g$. Since $d\varphi$ is a homeomorphism, there is $v\in \C^n$ such that $g_*\circ d\varphi(v)\not=0$. Thus, $d\phi\circ g_*\circ d\varphi(v)\not=0$. Therefore, $\lim\limits_{j\to +\infty}f(t_jv)$ is finite and non-zero, and thus ${\rm ord}_0f\leq k$, which shows that ${\rm ord}_0f= k$. In this case, $\lim\limits_{j\to +\infty}\frac{f(t_jv)}{t_j^k}=f_*(v)$ for all $v\in \C^n$. Therefore,
$f_*(v)=d\phi\circ g_*\circ d\varphi(v)$ for all $v\in \C^n$.

Since $f_*=d\phi\circ g_*\circ d\varphi$, $f_*^{-1}(1)$ and $g_*^{-1}(c)$ are homeomorphic, where $c=(d\phi)^{-1}(1)$. Since $g_*^{-1}(c)$ and $g_*^{-1}(1)$ are homeomorphic, we obtain, in particular, that $f_*^{-1}(1)$ and $g_*^{-1}(1)$ are homotopy equivalent.
\end{proof}

As a consequence of the above proposition and \cite[Theorem 3.2]{Sampaio:2016}, we obtain the following:

\begin{corollary}\label{cor:bi-Lip_A_equiv}
 Let $f,g\colon (\C^n,0)\to (\C,0)$ be two holomorphic function-germs. Assume that there are bi-Lipschitz homeomorphisms $\varphi\colon (\C^n,0)\to (\C^n,0)$ and $\phi\colon (\C,0)\to (\C,0)$ such that $f=\phi\circ g\circ \varphi$. Then the Milnor fibres of $f_*$ and $g_*$ are homotopy equivalent.
\end{corollary}

\subsection{Case of blow-spherical right equivalence}\label{subsec:bs}

Let us consider the {\bf spherical blowing-up} (at 0) of $\R^{n}$, $\rho \colon\mathbb{S}^n\times [0,+\infty )\to \R^{n}$), given by $\rho_{n}(x,t)=tx$.

Note that $\rho\colon\mathbb{S}^{n-1}\times (0,+\infty )\to \R^{n}\setminus \{0\}$ is a homeomorphism with inverse mapping $\rho_{n}^{-1}\colon\R^{n}\setminus \{0\}\to \mathbb{S}^{n-1}\times (0,+\infty )$ given by $\rho_{n}^{-1}(x)=(\frac{x}{\|x\|},\|x\|)$.

 The {\bf strict transform} of the subset $X$ under the spherical blowing-up $\rho_{n}$ is $X':=\overline{\rho_{n}^{-1}(X\setminus \{0\})}$. The subset $X'\cap (\mathbb{S}^{n-1}\times \{0\})$ is called the {\bf boundary} of $X'$ and is denoted by $\partial X'$.

\begin{definition}\label{def:bs_homeomorphism}
Let $X$ and $Y$ be subsets in $\mathbb{R}^n$ and $\mathbb{R}^m$ respectively. A homeomorphism $\varphi:(X,0)\rightarrow (Y,0)$ is said a {\bf blow-spherical homeomorphism} (at $0$), if  the homeomorphism 
$$\rho^{-1}_{m}\circ \varphi\circ \rho_{n} \colon X'_p\setminus \partial X'\rightarrow Y'\setminus \partial Y'$$
extends to a homeomorphism $\varphi'\colon X'\rightarrow Y'$. 
\end{definition}

\begin{proposition}\label{prop:bs_R-equiv}
 Let $f,g\colon (\C^n,0)\to (\C,0)$ be two reduced holomorphic function-germs. Let $\varphi\colon (\C^n,0)\to (\C^n,0)$ be a blow-spherical homeomorphism such that $f=g\circ \varphi$. Then the Milnor fibres of $f_*$ and $g_*$ are homotopy equivalent.
\end{proposition}
\begin{proof}
Let $k={\rm ord}_0 f$ and $m={\rm ord}_0 g$. Suppose by contradiction that $m>k$.

Since $\varphi$ is a blow-spherical homeomorphism, there is a homeomorphism $\nu_{\varphi}\colon \mathbb{S}^{2n-1}\to \mathbb{S}^{2n-1}$ given by $\nu_{\varphi}(x)=\lim\limits_{t\to 0^+}\frac{\varphi(tx)}{\|\varphi(tx)\|}$ and such that $\nu_{\varphi}(V(f_*)\cap \mathbb{S}^{2n-1})=V(g_*)\cap \mathbb{S}^{2n-1}$. So, the homeomorphism $d_0\varphi\colon \C^n\to \C^n$ given by
$$
d_0\varphi(x)=\left\{\begin{array}{ll}
   \|x\|\nu_{\varphi}(\frac{x}{\|x\|}),  &  \mbox{if }x\not =0\\
   0,  & \mbox{if }x=0
\end{array}\right.
$$
satisfies $d_0\varphi(V(f_*))=V(g_*)$.
Then, the function $\phi\colon \mathbb{S}^{2n-1}\setminus V(f_*)\to \C$ given $\phi(x)=\frac{f_*(x)}{g_*(\nu_{\varphi}(x))}$ is well-defined and continuous. Since 
$$
f_*(x)=\lim\limits_{t\to 0^+}\frac{f(tx)}{t^k}=\lim\limits_{t\to 0^+}\frac{g\circ \varphi(tv)}{t^k}=\lim\limits_{t\to 0^+}\frac{\|\varphi(tx)\|^m}{t^k}\frac{g\circ \varphi(tx)}{\|\varphi(tx)\|^m}
$$
and
$$
g_*(\nu_{\varphi}(x))=\lim\limits_{t\to 0^+}\frac{g\circ \varphi(tx)}{\|\varphi(tx)\|^m},
$$
we obtain that $\phi$ also satisfies the following
$$
\phi(x)=\lim\limits_{t\to 0^+}\frac{\|\varphi(tx)\|^m}{t^k}.
$$
In particular, the image of $\phi$ is contained in $(0,+\infty)$. Moreover, 
$$
f_*(x)=\textstyle{\|x\|^{k-m}\phi\big(\frac{x}{\|x\|}\big)g_*(d_0\varphi(x))=g_*\big(\|x\|^{\frac{k}{m}-1}\phi\big(\frac{x}{\|x\|}\big)^{\frac{1}{m}}d_0\varphi(x)\big)}
$$
for all $x\in \C^n\setminus V(g_*)$.

Similarly, the function $\tilde \phi\colon \mathbb{S}^{2n-1}\setminus V(g_*)\to \C$ given by $\tilde\phi(x)=\frac{g_*(x)}{f_*(\nu^{-1}_{\varphi}(x))}$ is well-defined, continuous, its image is contained in $(0,+\infty)$ and it satisfies
$$
g_*(x)=\textstyle{\|x\|^{m-k}\tilde\phi\big(\frac{x}{\|x\|}\big)f_*((d_0\varphi)^{-1}(x))=f_*\big(\|x\|^{\frac{m}{k}-1}\tilde\phi\big(\frac{x}{\|x\|}\big)^{\frac{1}{k}}(d_0\varphi)^{-1}(x)\big)}.
$$

Thus,
\begin{eqnarray*}
 f_*(x)&=&\textstyle{\phi\big(\frac{x}{\|x\|}\big)g_*(d_0\varphi(x))}\\
       &=& \textstyle{\phi\big(\frac{x}{\|x\|}\big)\tilde\phi\big(\frac{d_0\varphi(x)}{\|d_0\varphi(x)\|}\big)f_*(x)}.
\end{eqnarray*}
In particular, $\textstyle{\phi\big(\frac{x}{\|x\|}\big)\tilde\phi\big(\frac{d_0\varphi(x)}{\|d_0\varphi(x)\|}\big)}=1$ for all $x\in \C^n\setminus V(f_*)$.

Let $F=f^{-1}(1)$ and $G=g^{-1}(1)$ and consider $h\colon F\to G$ given by $h(x)=\|x\|^{\frac{k}{m}-1}\phi\big(\frac{x}{\|x\|}\big)^{\frac{1}{m}}d_0\varphi(x)$. We have that $h$ is a homeomorphism with inverse given by $h^{-1}(x)=\|x\|^{\frac{m}{k}-1}\tilde\phi\big(\frac{x}{\|x\|}\big)^{\frac{1}{k}}(d_0\varphi)^{-1}(x)$. 

Therefore, $f_*^{-1}(1)$ and $g_*^{-1}(1)$ are homeomorphic and, in particular, they are homotopy equivalent.
\end{proof}

\begin{corollary}\label{cor:bs_R-equiv}
 Let $f,g\colon (\C^n,0)\to (\C,0)$ be two reduced holomorphic function-germs. Let $\varphi\colon (\C^n,0)\to (\C^n,0)$ be a homeomorphism such that $f=g\circ \varphi$. Assume that $\tilde\varphi=\varphi\times id\colon (\C^n\times \C,0)\to (\C^n\times \C,0)$ is also a blow-spherical homeomorphism. Then $m(V(f),0)=m(V(g),0)$ and, moreover, the Milnor fibres of $f_*$ and $g_*$ are homotopy equivalent.
\end{corollary}
\begin{proof}
Let $k={\rm ord}_0 f$ and $m={\rm ord}_0 g$. Suppose by contradiction that $m>k$. In this case, $n>2$ and $k>1$.

Let us define $\tilde f, \tilde g\colon (\mathbb{C}^{n}\times\mathbb{C},0)\to (\mathbb{C},0)$ by
$\tilde f(x,t)=f(x)+t^k$ and $\tilde g(x,t)=g(x)+t^k$. Then $C(V(\tilde g),0)=\{t^k=0\}=\mathbb{C}^{n}\times\{0\}$, $\tilde f=\tilde g\circ\tilde \varphi$ and $C(V(\tilde f),0)$ is not a linear subspace. Since $\tilde \varphi$ is a blow-spherical homeomorphism, by \cite[Proposition 3.3]{Sampaio:2020b}, there is a homeomorphism $\psi\colon \C^{n+1}\to \C^{n+1}$ such that $\psi(C(V(\tilde f),0))=C(V(\tilde g),0)$. By \cite[Theorem 4.3]{Sampaio:2020b}, $C(V(\tilde f),0)$ is a linear subspace, which is a contradiction.

Therefore $m(V(f),0)=m(V(g),0)$.

Since $\tilde \varphi$ is a blow-spherical homeomorphism, $\varphi$ is also a blow-spherical homeomorphism. By Proposition \ref{prop:bs_R-equiv}, the Milnor fibres of $f_*$ and $g_*$ are homotopy equivalent.
\end{proof}

\section{Some results about Conjecture \ref{conj:weaker_problem_b}}

\subsection{Conjecture \ref{conj:weaker_problem_b} has a positive answer in dimension 2}

Using the bi-Lipschitz invariance of the tangent cones (see \cite{Sampaio:2016} and \cite{Sampaio:2015}) and the bi-Lipschitz classification of complex analytic curves, we directly obtain the following result.

\begin{proposition}\label{prop:conj_weaker_problem_b_dim_two}
 Let $f,g\colon (\C^2,0)\to (\C,0)$ be two reduced holomorphic function-germs. If there is a homeomorphism $\varphi\colon(\mathbb{C}^2,0)\to (\mathbb{C}^2,0)$ such that $f=g\circ \varphi$, then $f_*^{-1}(0)\setminus \{0\}$ and $g_*^{-1}(0)\setminus \{0\}$ have the same homology.
\end{proposition}

\subsection{Conjecture \ref{conj:weaker_problem_b} implies Conjecture \ref{conj:zariski_right}}

\begin{proposition}\label{prop:conj_weaker_problem_b_implies_zariski}
If Conjecture \ref{conj:weaker_problem_b} has a positive answer, then Conjecture \ref{conj:zariski_right} also has a positive answer.
\end{proposition}
\begin{proof}
Assume that Conjecture \ref{conj:weaker_problem_b} has a positive answer.

Assume by contradiction that there are two reduced holomorphic function-germs $f,g\colon(\mathbb{C}^n,0)\to (\mathbb{C},0)$ and a homeomorphism $\varphi\colon(\mathbb{C}^n,0)\to (\mathbb{C}^n,0)$ such that $f=g\circ \varphi$ and $k={\rm ord}_0 f<{\rm ord}_0 g=m$. Since Conjecture \ref{conj:zariski_right} has a positive answer in dimension 2, we have $n\geq 3$.

By A'Campo-Le's Theorem, we have $k>1$. 
Let $\bar f, \bar g\colon (\mathbb{C}^{n}\times\mathbb{C}^2,0)\to (\mathbb{C},0)$ given by
$\bar f(x,s_1,s_2)=f(x)+s_1^{m}+s_2^{m}$ and $\tilde g(x,s_1,s_2)=g(x)+s_1^{m}+s_2^{m}$. By changing $f$ and $g$ by $\bar f$ and $\bar g$ respectively, if necessary, we may assume that $n\geq 5$.

Let us define $\tilde f, \tilde g\colon (\mathbb{C}^{n}\times\mathbb{C},0)\to (\mathbb{C},0)$ by
$\tilde f(x,t)=f(x)+t^k$ and $\tilde g(x,t)=g(x)+t^k$. Then $C(V(\tilde g),0)=\{t^k=0\}=\mathbb{C}^{n}\times\{0\}$ and $\tilde f=\tilde g\circ\tilde \varphi$. In particular, $H^i(C(V(\tilde g),0)\setminus \{0\})=0$ for any $i\in \{1,...,2n-2\}$.

Since we are assuming that Conjecture \ref{conj:weaker_problem_b} has a positive answer, we have $H^i(C(V(\tilde f),0)\setminus \{0\})=0$ for any $i\in \{1,...,2n-2\}$. 
\begin{claim}\label{claim:linear_subspace}
$C(V(\tilde f),0)$ is a linear subspace.
\end{claim}
\begin{proof}[Proof of Claim \ref{claim:linear_subspace}]
This claim holds true in a much more general setting; see \cite[Claim 3.1.3]{Sampaio:2024}. For completeness, we present a proof here. 

Let $\pi\colon V^*\to\widetilde{V}$ denote the quotient map by the $\mathbb{C}^*$-action, where $V=C(V(\tilde f),0)$ and $V^*=V\setminus\{0\}$.
Let $\mathbb{P}=\C P^{n}$ and $U=\mathbb{P}\setminus \widetilde{V}$.
 By the Alexander Duality, 
 $$
 H^s(\mathbb{P},\widetilde{V})=H_{2n-s}(U).
 $$

Since $U$ is a smooth affine hypersurface of complex dimension $n$, it follows from the Andreotti-Frankel theorem (see \cite{AndreottiF:1959}) that $H_m(U)=0$ for all $m\geq n+1$.

Since $n\geq 5$, $H^2(\mathbb{P},\widetilde{V})=H^3(\mathbb{P},\widetilde{V})=0$, and thus the morphism $j^2\colon H^2(\mathbb{P})\to H^2(\widetilde{V})$, induced by the inclusion $j\colon \widetilde{V}\to \mathbb{P}$, is an isomorphism.

We know that the integral cohomology algebra $H^{\bullet}(\mathbb{P})$ is a truncated polynomial algebra $H^{\bullet}(\mathbb{P})=\mathbb{Z}[\alpha]/(\alpha^{N})$ generated by an element $\alpha$ of degree 2 (see \cite[Chapter 5, Proposition 1.6]{Dimca:1992}).
Thus, $\alpha_V=j^2(\alpha)$ is a generator of $H^2(\widetilde{V})$.

The associated Gysin sequence in cohomology leads to the following 
\[\begin{tikzcd}
   ... \rightarrow H^{2p+1}(X_0)\rightarrow  H^{2p}(\widetilde{V})\ar{r}{\psi_p}&H^{2p+2}(\widetilde{V})\rightarrow  H^{2p+2}(X_0)\rightarrow ...
\end{tikzcd}\]
where $\psi_p$ is the cup product with $\alpha_V$. Since $H^i(X_0)=0$ for all $1\leq i\leq 2n-2$, we obtain that $\psi_p$ is an isomorphism for all $p\in\{1,...,n-2\}$ and, in particular, $H^{2n-2}(\widetilde{V})=\mathbb{Z}$ and $\alpha_V^{n-1}$ is a generator of $H^{2n-2}(\widetilde{V})$. This also shows that $\widetilde{V}$ is an irreducible algebraic set.

Let $E_m$ be a generic hyperplane in $\mathbb{P}$ of dimension $m$. Then ${\rm deg}(\widetilde{V})=\#(\widetilde{V}\cap E_{N-d})$ and $[\widetilde{V}]={\rm deg}(\widetilde{V})[E_{n-1}]$ in $H^{2(n-1)}(\mathbb{P})$, where $[A]$ denotes the fundamental class of $A$. Hence,
$$
\langle j^{2(n-1)}(\alpha^{n-1}),[\widetilde{V}]\rangle=\langle \alpha_V^{n-1},[\widetilde{V}]\rangle=\langle \alpha_V^{n-1},{\rm deg}(\widetilde{V})[E_{n-1}]\rangle={\rm deg}(\widetilde{V}).
$$
Therefore $\alpha_V^{n-1}={\rm deg}(\widetilde{V}) \cdot g$, where $g$ is a generator of $H^{2n-2}(\widetilde{V})$. Since $\alpha_V^{n-1}$ is also a generator of $H^{2n-2}(\widetilde{V})$, we obtain ${\rm deg}(\widetilde{V})=1$, and therefore $V$ is a linear subspace.

\end{proof}

However, since $f_*\not \equiv 0$ and $k>1$, we obtain that $f_*+t^k$ cannot be a power of a linear form and, in particular, $C(V(\tilde f),0)=V(f_*+t^k)$ is not a linear subspace, which is a contradiction.
\end{proof}

\section{Zariski Multiplicity Conjecture mod 2}

\begin{proposition}\label{prop:equiv_zariski}
If Conjecture Z mod 2 has a positive answer, then Conjecture \ref{conj:zariski_right} also has a positive answer.
\end{proposition}
\begin{proof}
Suppose that Conjecture Z mod 2 has a positive answer.

Assume by contradiction that Conjecture \ref{conj:zariski_right} does not have a positive answer. Thus, there are two reduced holomorphic function-germs $f,g\colon(\mathbb{C}^n,0)\to (\mathbb{C},0)$ and a homeomorphism $\varphi\colon(\mathbb{C}^n,0)\to (\mathbb{C}^n,0)$ such that $f=g\circ \varphi$, but $m(V(f),0)\not=m((V(g),0)$. Since Conjecture \ref{conj:zariski} has a positive answer in dimension $2$, then $n\geq 3$.

We assume that $m:=m(V(f),0)>k:=m((V(g),0)$, the case $m(V(f),0)<m((V(g),0)$ is proved in the same way. Thus, we define $\tilde f,\tilde g\colon (\mathbb{C}^n\times \C,0)\to (\mathbb{C},0)$ given by $\tilde f(x,t)=f(x)+t^{k+1}$ and $\tilde g(x,t)=g(x)+t^{k+1}$.  Then, the homeomorphism $\tilde\varphi\colon(\mathbb{C}^n\times \C,0)\to (\mathbb{C}^n\times \C,0)$ given by $\tilde \varphi(x,t)=(\varphi(x),t)$ satisfies $\tilde f=\tilde g\circ \tilde \varphi$. Since we are assuming that Conjecture Z mod 2 has a positive answer, then $m(V(\tilde f),0)\equiv m(V(\tilde g),0)\, {\rm mod}\, 2$. But $m(V(\tilde f),0)=k+1$ and $m(V(\tilde g),0)=k$, which is a contradiction. Therefore, Conjecture \ref{conj:zariski_right} has a positive answer.
\end{proof}

% \noindent {\bf Acknowledgements.}

\end{document}